\documentclass{amsart}
\usepackage[english]{babel}
\usepackage[cp1250]{inputenc}
\usepackage{amssymb,amsthm,amsfonts}
\usepackage{verbatim}
\usepackage{url}

\newcommand{\RR}{\mathbb R}
\newcommand{\FF}{\mathbb F}

\newcommand{\NN}{\mathbb N}

\newcommand{\CC}{\mathbb C}

\newcommand{\cB}{\mathcal B}
\newcommand{\cI}{\mathcal I}

\newcommand{\cK}{\mathcal K}

\newcommand{\cH}{\mathcal H}

\newcommand{\cS}{\mathcal S}

\newcommand{\cA}{\mathcal A}

\newcommand{\cY}{\mathcal Y}
\newcommand{\cZ}{\mathcal Z}

\newcommand{\im}{\mathrm{Im}}
\newcommand{\re}{\mathrm{Re}}
\newcommand{\benu}{\begin{enumerate}}
\newcommand{\eenu}{\end{enumerate}}
\newcommand{\bop}{\begin{opomba}}
\newcommand{\eop}{\end{opomba}}

\newcommand{\Bor}{\mathrm{Bor}}

\newcommand{\B}{\mathrm{Bor}}

\newcommand{\id}{\mathrm{Id}}

\newcommand{\beqn}{\begin{eqnarray*}}
\newcommand{\eeqn}{\end{eqnarray*}}

\newtheorem{corollary}{Corollary}
\newtheorem{theorem}{Theorem}
\newtheorem{proposition}{Proposition}
\newtheorem{lemma}{Lemma}
\theoremstyle{definition}
\newtheorem{definition}{Definition}

\newtheorem{remark}{Remark}

\theoremstyle{remark}

\newcommand{\bdefi}{\begin{definition}}
\newcommand{\edefi}{\end{definition}}
\newcommand{\bcor}{\begin{corollary}}
\newcommand{\ecor}{\end{corollary}}
\newcommand{\bthe}{\begin{theorem}}
\newcommand{\ethe}{\end{theorem}}
\newcommand{\bpro}{\begin{proposition}}
\newcommand{\epro}{\end{proposition}}
\newcommand{\blem}{\begin{lemma}}
\newcommand{\elem}{\end{lemma}}
\newcommand{\brem}{\begin{remark}}
\newcommand{\erem}{\end{remark}}
\newcommand{\bequ}{\begin{equation}}
\newcommand{\eequ}{\end{equation}}
\newcommand{\bprf}{\begin{proof}}
\newcommand{\eprf}{\end{proof}}

\begin{document}
\title{Non-negative spectral measures and representations of $C^\ast$-algebras}

\author{Alja\v z Zalar}

\address{Alja\v z Zalar, University of Ljubljana, Faculty of Math.~and Phys., Dept.~of Math.}
\email{aljaz.zalar@imfm.si}

\date{\today}

\begin{abstract} 
	Regular normalized $W$-valued spectral measures on a compact Hausdorff space $X$ are in one-to-one correspondence with 
	unital $\ast$-rep\-re\-sen\-ta\-tions 	$\rho:C(X,\CC)\to W$, where $W$ stands for a von Neumann algebra.
	In this paper we show that for every compact Hausdorff space $X$ and every von Neumann algebras $W_1,W_2$
	there is a one-to-one correspondence between unital $\ast$-representations 	
		$\rho:C(X,W_1)\to W_2$ and special $B(W_1,W_2)$-valued measures on $X$ that we call
	non-negative spectral measures.
	Such measures are special cases of non-negative measures that we introduced in our previous paper \cite{Cim-Zal}
	in connection with moment problems for operator polynomials.
\end{abstract}



\maketitle


\section{Introduction}

A $\ast$-representation of a $C^\ast$-algebra $\cA$ is an algebra homomorphism $\rho:\cA\to W$ such
that $\rho(a^{\ast})=\rho(a)^{\ast}$ for every $a\in A$, where $W$ is a von Neumann algebra.
Our main result is the following theorem on $\ast$-representations of the form
$\rho: C(X,W_1) \to W_2$, where $X$ is a compact Hausdorff space and $W_1$, $W_2$ are von Neumann algebras.
It is a generalization of the usual situation, i.e., $\ast$-representations of the form $\rho: C(X,\CC)\to W$ (see Theorem \ref{reprezentacija} below). 
By $B(W_1,W_2)$ we denote the Banach space of all bounded linear operators from $W_1$ to $W_2$.

\bthe \label{reprezentacijski-izrek-0}
	Let $X$, $W_1, W_2$, $B(W_1,W_2)$ be as above and $\rho: C(X,W_1) \to W_2$ a linear map.
	Let $\B(X)$ be a Borel $\sigma$-algebra on $X$.
	The following statements are equivalent.
	\benu
		\item $\rho: C(X,W_1) \to W_2$ is a unital $\ast$-representation. 
		\item	There exists a unique regular normalized non-negative spectral measure $M:\B(X)\to B(W_1,W_2)$ such that
				$$\rho(F)=\int_X F\; dM$$
			for every $F\in C(X,W_1)$.
	\eenu
\ethe

A set function  
		$$M \colon \B(X) \to B(W_1,W_2)$$ 
	is \textit{a non-negative spectral measure} if for every hermitian projection $P\in W_1$ the set function
		$$M_P \colon \B(X) \to W_2, \quad M_P(\Delta):=M(\Delta)(P)$$
	is a spectral measure such that 
	$$M_P(\Delta_1)M_Q(\Delta_2)=M_{PQ}(\Delta_1\cap \Delta_2) $$
holds for all hermitian projections $P,Q\in W_1$ and all sets $\Delta_1, \Delta_2\in \B(X)$.

\brem
	\benu
	\item
		Spectral measures and their adaptations are well-studied in the representation theory (e.g., \cite{Bas-Fer-Kar}, \cite{Mas}, \cite{Mas-Ros}, \cite{Ros}).
		We introduced non-negative measures in \cite{Cim-Zal}, where we studied moment problems in the case of operator polynomials (see Sections 3-5 below for 	
		a concise treatment of non-negative measures).
		Non-negative spectral measures are their special cases (see Sections 7, 8).
	\item 
		Note that since every von Neumann algebra is a dual of a Banach space, the existence of a representing measure in Theorem \ref{reprezentacijski-izrek-0} is already covered as a special case of 	
		\cite[Theorem 3.3.]{Mitt-Young}.
		The interesting part of Theorem \ref{reprezentacijski-izrek-0} is a concrete description of the representing measure in this special case and
		a one-to-one correspondence between $\ast$-representations and measures.
	\eenu
\erem

The paper is structured in the following way. In Section 2 we introduce some terminology and state a well-known representation theorem for abelian $C^{\ast}$-algebras. In Section 3 we present the complex version of the measure and integration theory from \cite{Cim-Zal} in a more systematic way.
Section 4 provides a characterization of non-negative measures (see Theorem \ref{karakterizacija-nenegativnih-kompleksnih-mer}).
In Section 5 we extend the integration theory to a Banach space which in particular constists of all bounded measurable $W_1$-valued functions 
and obtain a slight extension of \cite[Proposition 2]{Cim-Zal}; see Theorem \ref{tretji-Rieszov-operatorski-izrek}.
In Section 6 we show how our measures are connected with the measures from \cite{Mitt-Young} (see Proposition \ref{opomba-Mitter-Young}).
In Section 7 we introduce non-negative spectral measures. Section 8 provides a characterization of non-negative spectral measures (see Theorem \ref{karakterizacija-nenegativnih-spektralnih-mer}),
which is then used in Section 9 to prove Theorem \ref{reprezentacijski-izrek-0} (see Theorem \ref{reprezentacijski-izrek-1} and Corollary \ref{reprezentacijski-izrek-2}).


\section{ Preliminaries }
\label{preliminaries} 
 
Let $(X,\cS,\cH)$ be a measure space, i.e., $X$ is a set, $\cS$ a $\sigma$-algebra on $X$
and $\cH$ a Hilbert space, and $\id_{\cH}$ denotes the identity operator on $\cH$.
\textsl{Spectral measure} $F:\cS\to B(\cH)$ is a positive operator-valued measure with an additional property that it maps
into the set of hermitian projections; see \cite[Definition 2]{Ber}. $F$ is
\textsl{normalized} if $F(X)=\id_{\cH}$. $F$ on a locally compact space $X$, equipped with a Borel $\sigma$-algebra $\B(X)$, is \textsl{regular} if the complex measures 
	$$F_{h_1,h_2}: \B(X)\to \CC,\quad F_{h_1,h_2}(\Delta):=\left\langle F(\Delta)h_1,h_2\right\rangle$$ 
are regular for all $h_1,h_2\in\cH$. It is a well-known fact that every normal operator $A$ can be represented as an integral with respect to a unique regular normalized spectral measure $F$, i.e., $A=\int_\RR t\; dF(t)$. 

A $\ast$-representation of a $C^\ast$-algebra $\cA$ is an algebra homomorphism $\rho:\cA\to W$ such
that $\rho(a^{\ast})=\rho(a)^{\ast}$ for every $a\in A$, where $W$ is a von Neumann algebra. Spectral measures are interesting also due to the following result; see \cite[p.~259]{Con} and note that $B(\cH)$ can be replaced by $W$ by \cite[Theorem 2.7.4]{Pet}.

\bthe \label{reprezentacija}
	Let $X$ be a compact Hausdorff space, $W$ a von Neumann algebra and $\rho: C(X,\CC)\to W$ a linear map. 
	Let $\B(X)$ be a Borel $\sigma$-algebra on $X$.
	The following statements are equivalent.
	\benu
		\item $\rho: C(X,\CC)\to W$ is a unital $\ast$-representation. 
		\item There exists a unique regular normalized spectral measure $F:\Bor(X)\to W$ such that
				$\rho(f)=\int_X f\; dF$
			for every $f\in C(X,\CC)$. 
	\eenu
\ethe


\brem \label{regularnost}
	\benu
		\item The assumptions that $X$ \textsl{is a compact Hausdorff space} and $\rho$ \textsl{a linear map of the form} $\rho: C(X,\CC)\to 
			W$ can be replaced by the assumptions that $X$ \textsl{is a locally compact Hausdorff space} and
			$\rho$ \textsl{a linear map of the form} $\rho: C_0(X,\CC)\to W$, where $C_0(X,\CC)$ denotes the space of functions vanishing 
			at infinity. By compactifying $X$ with one point to $X_\infty$ and using Theorem \ref{reprezentacija}, 
			$\ast$-representations of the form $\rho: C_0(X,\CC)\to W$ are in one-to-one correspondence with
			the regular normalized spectral measures $F:\B(X_{\infty})\to W$, where $\B(X_{\infty})$ is the Borel $\sigma$-algebra on $X_{\infty}$. 
			However, this result is also covered by \cite[Theorem 4.1.]{Mitt-Young}.
		\item The Baire $\sigma$-algebra is a $\sigma$-algebra generated by all compact subsets of $X$, which are $G_{\delta}$ sets, i.e., a
			countable intersection of open sets. \cite[Theorem 19]{Ber} is the same result as Theorem \ref{reprezentacija}, where the Borel
			$\sigma$-algebra is replaced by the Baire $\sigma$-algebra. In general one has to be cautious when working with Baire or 
			Borel $\sigma$-algebras. For $\sigma$-compact and metrizable spaces they coincide, but for general topological spaces this is not the
			case. The reason for Borel $\sigma$-algebra being appropriate in Theorem \ref{reprezentacija} is the following lemma (see \cite[Proposition V.4.1]{Con}) 
			and working with nets instead of sequences.
			\blem \label{density}
				The ball in $C(X,\CC)$ is a dense subset of the ball in $C(X,\CC)^{\ast\ast}$ equipped with a $weak^{\ast}$-topology.
			\elem
	\eenu
\erem


\section{Non-negative measures} \label{razdelek-z-nenegativnimi-merami}

For Banach spaces $\cY$, $\cZ$ we denote by $B(\cY,\cZ)$ the Banach space of all bounded linear operators from $\cY$ to $\cZ$.
In the case $\cY=\cZ$ we write $B(\cY)$ for $B(\cY,\cY)$. 
For a von Neumann algebra $W\subseteq B(\cH)$, where $\cH$ is a Hilbert space, we denote by $W_h$, $W_+$ the subsets of $W$ of all hermitian operators and all positive operators respectively. By a positive operator we mean a hermitian operator $A$, which satisfies $\left\langle Ah,h \right\rangle\geq 0$ for every $h\in\cH$ (Here $\left\langle \cdot,\cdot\right\rangle$ denotes the inner product on $\cH$.).

	Let $X$ be a set, $\cS$ a $\sigma$-algebra on $X$, $\cH$, $\cK$ Hilbert spaces over $\FF\in\left\{\RR, \CC\right\}$ and $W_1\subseteq B(\cH), W_2\subseteq B(\cK)$ von Neumann algebras.
	For $\FF=\CC$ ($\FF=\RR$) a set function 
		$$m \colon \cS \to B(W_1,W_2) \quad \left(m \colon \cS \to B\left((W_1)_h,(W_2)_h\right)\right)$$
	is \textit{a non-negative measure} if for every $A\in (W_1)_+$ the set function 
		$$m_A \colon \cS \to W_2, \quad m_A(\Delta):=m(\Delta)(A),$$
	is a positive operator-valued measure.
A quadruple $(X,\cS,W_1,W_2)$ is \textsl{a measure space} and 
a pentuple $(X,\cS,W_1,W_2,m)$ \textsl{a space with a measure $m$}.

\brem
	The reason for the distinction in the definition of $m$ between the real and the complex case lies in the fact, that in the complex case every element
	$A\in W_1$ can be written as a $\CC$-linear combination of two hermitian elements, i.e., $A=\frac{A+A^\ast}{2}+i\frac{A-A^{\ast}}{2i}$, while this is not true
	in the real case. Hence, in the real case for the uniqueness of $m$ it is not sufficient to know all set functions $m_A$ for every $A\in (W_1)_+$. However, it 	
	suffices if $m$ is of the form $m \colon \cS \to B\left((W_1)_h,(W_2)_h\right)$.
\erem

Let $(X,\cS,W_1,W_2,m)$ be a space with a measure $m$.
A $\cS$-measurable complex function $f:X\to\CC$ is \textsl{$m$-integrable}, if it is $m_A$-integrable for every $A \in (W_1)_+$.
The set of all \textsl{$m$-integrable} functions is denoted by $\cI(m)$.

\brem \label{pred-konvergencnim-izrekom}
			Given a positive operator-valued measure $E:\cS\to B(\cK)$, where $\cK$ is a Hilbert space, a $\cS$-measurable function $f:X\to \CC$ is called 
			$E$-integrable, if 
			there exists a constant $K_f\in \RR$ such that $\int_{X} \left|f\right|\;dE_k\leq K_f \left\|k\right\|^2$ for every $k\in \cK$.
			Here $E_k$ denotes a positive measure $E_k:\cS\to [0,\infty)$ defined by $E_k(\Delta):=\left\langle E(\Delta)k,k\right\rangle$ for every $\Delta\in\cS$.
			Then in the case $\FF=\RR$ the mapping $(k_1,k_2)\mapsto \frac{1}{4}\left(\int_X f\; dE_{k_1+k_2}-\int_X f\;dE_{k_1-k_2}\right)$ is a bounded bilinear form, 
			while in the case 
			$\FF=\CC$ the mapping $(k_1,k_2)\mapsto \frac{1}{4}\sum_{j=0}^{3}i^j\int_X f\;dE_{k_1+i^j k_2}$ is a bounded sesquilinear form.
\erem

The set $\cI(m)$ is a complex vector space and it consists of at least all bounded $\cS$-measurable complex functions.
In particular, for $\cS=\B(X)$ we have $C_c(X,\CC)\subset \cI(m)$.

The following convergence theorem will be frequently used in the sequel.

\bthe \label{konvergencni-izrek}
	Let $\{f_n\}_{n\in\NN}$ be an increasing sequence of positive $E$-integrable functions that pointwise 
	converges to a $\cS$-measurable function $f$. 
	If there exists $B\in B(\cK)$ such that $\int_X f_n\;dE\preceq B$,	
	then $f$ is $E$-integrable and 
		$$\lim_n \int_X f_n\; dE = \int_X f\;dE,$$ 
	where the limit is taken in the strong operator topology.
\ethe

\bprf
	Since by the usual convergence theorem we have 
		$\int_X f\,dE_k=$$\lim_{n}\int_X f_n \;dE_k $$\leq \left\langle Bk,k \right\rangle$
	for every $k\in\cK$, $f$ is $E$-integrable (take $K_f=\left\|B\right\|$ in Remark \ref{pred-konvergencnim-izrekom}).
	Then proceed as in the proof of \cite[Theorem 11(iii)]{Ber}.
\eprf

Given $A\in W_1$ we write $\re(A):=\frac{1}{2}(A+A^{\ast})\in W_1$ and $\im(A):=\frac{i}{2}(A^{\ast}-A)\in W_1$ for its the real and imaginary part, while for
$A\in (W_1)_h$ we write $A_+$ and $A_-$ for its positive and negative part ($A_+,A_-\in W$ by \cite[Proposition 2 on p.~3]{Dix}).

	For each $m$-integrable function $f$ and each operator $A\in W_1$ we define $\int_X {f\; dm_A}$ as
		$$\int_X {f\; dm_{\re(A)_+}}-\int_X {f\; dm_{\re(A)_-}}+ \\
			 i\cdot\int_X {f\; dm_{\im(A)_+}}-i\cdot\int_X {f\; dm_{\im(A)_-}}.$$
Let $\cI(m) \otimes_\FF W_1$ be an algebraic tensor product of $\cI(m)$ and $W_1$ over $\FF\in\left\{\RR,\CC\right\}$.
We define the map
	$$\cB:\cI(m) \times W_1 \to W_2, \quad \cB(f,A) = \int_X {f\; dm_A}.$$ 
Let $\cI(m)_+$ be the set of all functions $f\in\cI(m)$, such that $f(x)\geq 0$ for every $x\in X$. 	

\bpro \label{bilinearnost-preslikave-B}
	The map $\cB$ is bilinear.
\epro

\bprf
	It suffices to consider $f,g\in\cI(m)_+$, $A,B\in (W_1)_h$.
	Equality $\cB(\alpha f+ \beta g,A)=\alpha\cB(f,A) + \beta \cB(g,A)$ easily follows by the definitions.
	Equality $\cB(f,A+B)=\alpha\cB(f,A) + \beta \cB(f,B)$ is equivalent to the equality of
	$C:=\int_X f\; dm_{(A+B)_+} + \int_X f\; dm_{A_-}+ \int_X f\; dm_{B_-}$ and 
	$D:=\int_X f\; dm_{(A+B)_-} + \int_X f\; dm_{A_+}+ \int_X f\; dm_{B_+}$.
	There is an increasing sequence $\{s_k\}_{k\in\NN}$ of simple functions 
	$s_k\in \cI(m)_+$ such that $\lim_{k} s_k=f$. By Theorem \ref{konvergencni-izrek}, 
	$C=\lim_k\int_X s_k\; dm_{(A+B)_+} + \int_X s_k\; dm_{A_-}+ \int_X s_k\; dm_{B_-}=
	\lim_k\int_X s_k\;  dm_{(A+B)_++A_-+B_-}=\int_X f\;  dm_{(A+B)_++A_-+B_-}.$
	Similarly $D=\int_X f\;  dm_{(A+B)_-+A_++B_+}$.
	By $(A+B)_++A_-+B_-=(A+B)_-+A_++B_+$, it follows $C=D$, which concludes the proof.
\eprf

By the universal property of the tensor product the bilinear map $\cB$ can be extended to the linear map
	$$	\bar \cB: \cI(m) \otimes_\FF W_1 \to W_2,\quad
			\bar \cB \left( F:=\sum_{i=1}^{n}f_i\otimes A_i\right) =  \sum_{i=1}^{n}  \int_X {f_i\; dm_{A_i}}=:\int_X F\; dm.$$

We call $F\in\cI(m)\otimes_\FF W_1$ \textsl{positive} if $F(x)\succeq 0$ for every $x\in X$ and write $F\succeq 0$.
In the following proposition we list some properties of the integral with respect to $m$.

\bpro\label{lastnosti-kompleksnega-integrala}
	Let $(X,\cS,W_1\subseteq B(\cH),W_2\subseteq B(\cK), m)$ be a space with a measure $m$ and $\FF\in\left\{\RR,\CC\right\}$. For all $F,G\in \cI(m)\otimes_\FF 
	W_1$, all operators $A\in W_1$, all numbers $\lambda\in\FF$ and all sets $\Delta\in\cS$ 
	the following equalities hold.
		\bequ \label{enakost1} \int_X (F+G)\;dm=\int_X F \;dm + \int_X G \;dm,\eequ
		\bequ \label{enakost2}\int_X \lambda F\;dm=\lambda \int_X F\;dm.\eequ
		\bequ \label{enakost3} \int_X \left(\chi_\Delta\otimes A\right)\;dm=m_A(\Delta).\eequ
	If $F\in \cI(m)\otimes_\FF W_1$ satisfies $F\succeq 0$, then
		\bequ \label{neenakost1}\int_X F\;dm\succeq 0.\eequ
\epro

\bprf
	(\ref{enakost1}), (\ref{enakost2}) follow by the construction of the map $\bar \cB$.
	It suffices to prove (\ref{enakost3}) for $A\in (W_1)_+$. Since $m_A$ is a positive operator-valued measure it follows 
		$\int_X \left(\chi_\Delta\otimes A\right)\;dm$$ = \int_X \chi_\Delta\; dm_A=m_A(\Delta).$
	It remains to prove (\ref{neenakost1}). Every $F\in \cI(m)\otimes_\FF W_1$, $F\succeq 0$ can be expressed as
		$\sum_{i=1}^{n_1} r_i \otimes B_i - \sum_{j=1}^{n_2} s_j \otimes C_j,$ 
	where $r_i \otimes B_i, s_j \otimes C_j \in \cI(m)_+ \otimes_\FF (W_1)_+$, $n_1,n_2\in\NN$. For 
	every $\ell\in \NN$ we define the set
	$X_{\ell}:=\left(\bigcap_i r_i^{-1}[0,\ell]\right) \bigcap \left(\bigcap_i s_j^{-1}[0,\ell]\right).$
	The sequence $X_{\ell}$ is increasing and $X=\cup_{\ell\in\NN} X_{\ell}$. For every $i,j$ there are positive simple
	functions $t^{1}_{\ell i}$, $t^{2}_{\ell j}$ such that
		$\left\|t^{1}_{\ell i} - \chi_{\ell} r_i \right\|_{\infty}\leq \frac{1}{2n_1 \ell \left\|B_i\right\|}$,
		$\left\|t^{2}_{\ell j}  -  \chi_{\ell} s_j \right\|_{\infty}\leq \frac{1}{2n_2 \ell \left\|C_j\right\|}$,
	where $\chi_\ell$ is a characteristic function of $X_{\ell}$.
	For every $\ell\in \NN$ we define 
	$G_{\ell}(x):=\left(\sum_{i=1}^{n_1} t^{1}_{\ell i} \otimes B_i - \sum_{j=1}^{n_2} t^{2}_{\ell j} \otimes C_j\right)$.
	Therefore 
		$\left\|\chi_{\ell} F - G_\ell\right\| \leq \frac{1}{\ell}$.
	Together with $\chi_{\ell} F\succeq 0$ it follows that
		$G_{\ell}\succeq -\frac{1}{\ell} \id_{\cH},$
	where $\id_\cH$ denotes the identity operator on $\cH$. 
	Each $G_{\ell}$ is of the form $\sum_{k} \chi_{\Delta_{k\ell}} \otimes D_{k\ell}$, where
	$\Delta_{k\ell}\in \cS$, $\Delta_{k\ell}\cap \Delta_{k'\ell}=\emptyset$ for $k\neq k'$, $\cup_k \Delta_{k\ell}=X$
	and $D_{{k\ell}}\succeq -\frac{1}{\ell} \id_{\cH}$.
	It follows that
 	\beqn 
		 \int_X {G_{\ell}\;dm}&=& \int_X {\left(\sum_{k} \chi_{\Delta_{k\ell}} \otimes D_{k\ell}\right)\;dm}= 
		 					\sum_{k} m_{D_{k\ell}}(\Delta_{k\ell})
		 \succeq \sum_k m_{-\frac{1}{\ell} \id_{\cH}}(\Delta_{k\ell})\\ &=&
							-\frac{1}{\ell} \sum_k m_{\id_{\cH}}(\Delta_{k\ell})
			=	-\frac{1}{\ell}\cdot  m_{\id_{\cH}}\left(\cup_k \Delta_{k\ell}\right)
			= -\frac{1}{\ell}\cdot m_{\id_{\cH}}(X)
	\eeqn
	Since for every $i,j$, the functions $r_i$, $s_j$ are positive and the sequence $X_\ell$ increases, the sequences $t^{1}_{\ell i}$, 	
	$t^{2}_{\ell j}$
	can be chosen such that they increase, i.e., for fixed $i,j$ we have $t^1_{1i}\leq t^1_{2i}\leq t^1_{3i}\leq \ldots$ and
	$t^2_{1j}\leq t^2_{2j}\leq t^2_{3j}\leq \ldots$. 
	By Theorem \ref{konvergencni-izrek},
		$\lim_{\ell}\int_{X} t^{1}_{\ell i}\; dm_{B_i} = \int_{X} r_i\;dm_{B_i},$
		$\lim_{\ell}\int_{X} t^{2}_{\ell j}\; dm_{C_j} = \int_{X} s_j\;dm_{C_j}.$
	It follows that
		$\int_{X} F\; dm = \lim_{\ell} \int_{X} G_{\ell}\; dm  \succeq \lim_{\ell} \left(-\frac{1}{\ell}\cdot m_{\id_{\cH}}(X)\right)=0.$
	This proves $\int_{X} F\; dm\succeq 0$, which is (\ref{neenakost1}).
\eprf


\section{Characterization of non-negative measures}

Let $(X,\cS,\cH)$ be a measure space, i.e., $X$ is a set, $\cS$ a $\sigma$-algebra on $X$
and $\cH$ a Hilbert space. \cite[Theorem 2]{Ber} characterizes positive operator-valued measures on $(X,\cS,\cH)$ via families $\left\{\mu_h\right\}_{h\in\cH}$ of finite positive measures.
We would like to have an analoguous characterization in the case of non-negative measures on a measure space $(X,\cS,W_1,W_2)$, i.e., $X$ is a set, $\cS$ a $\sigma$-algebra on $X$ and $W_1$, $W_2$ are von Neumann algebras.

The following theorem provides a characterization of non-negative measures on $(X,\cS,W_1,W_2)$ via families $\left\{E_A\right\}_{A\in (W_1)_+}$ of positive operator-valued measures. 

\bthe \label{karakterizacija-nenegativnih-kompleksnih-mer}
	Let $(X,\cS,W_1\subseteq B(\cH),W_2\subseteq B(\cK))$ be a measure space and $$\left\{E_A\right\}_{A\in (W_1)_+}$$ a family of 
	positive operator-valued measures $E_A:\cS\to W_2$.
	
	There exists a unique non-negative measure $m$ such that
		$$m_A=E_A$$ 
	for all operators $A\in (W_1)_+$
	iff the following conditions hold.
		\bequ \label{pogoj111}  E_{A+B}(\Delta)=E_A(\Delta)+E_B(\Delta),\eequ
 		\bequ \label{pogoj222} E_{\lambda A}(\Delta)=\lambda E_A(\Delta),\eequ
	for all operators $A,B\in (W_1)_+$, all real numbers $\lambda\in\RR^+$, and all sets $\Delta\in \cS$, and for each set $\Delta\in\cS$ there
	exists a constant $k_\Delta\in\RR^{>0}$ such that
		\bequ \label{pogoj333} \left\|E_A(\Delta)\right\|\leq k_\Delta \left\|A\right\|\eequ
	for all operators $A\in (W_1)_+$.
	
	Every family $\left\{E_A\right\}_{A\in (W_1)_+}$ which satisfies the conditions above is called 
	\textsl{a compatible family of positive operator-valued measures}.
\ethe

\bprf
	The nontrivial direction is the if part.
	We have to prove the well-definedness of the set function
		$$m \colon \cS \to B(W_1,W_2),$$
		$$m(\Delta)(A):= \left(E_{\re(A)_+}(\Delta) - E_{\re(A)_-}(\Delta)\right) + i \left(E_{\im(A)_+}(\Delta) - E_{\im(A)_-}(\Delta)\right),$$
	where $B(W_1,W_2)$ denotes the Banach space of all bounded linear operators from $W_1$ to $W_2$.
	For the well-definedness we have to show that for each $\Delta\in \cS$ the map $m^\Delta: W_1 \to W_2$, $m^\Delta(A):=m(\Delta)(A)$
	is linear and bounded. If Hilbert spaces $\cH$, $\cK$ are complex, then by the usual decompositions of $\lambda\in \CC$ and $A\in W_1$ into the real and imaginary 
	part it suffices to prove the $\RR$-linearity and the boundedness of $m^\Delta$ over $(W_1)_h$.
	
	\textsl{Additivity of $m^\Delta$.} For $A, B\in (W_1)_h$ the equality $m^\Delta(A+B)=m^\Delta(A)+m^{S}(B)$ is equivalent to 
		$$E_{(A+B)_+}(\Delta)-E_{(A+B)_-}(\Delta)=\left(E_{A_+}(\Delta)-E_{A_-}(\Delta)\right) + \left(E_{B_+}(\Delta)-E_{B_-}\right)(\Delta),$$ 
	which is further equivalent to
		$$E_{(A+B)_+}(\Delta)+ E_{A_-}(\Delta)+ E_{B_-}(\Delta) = E_{(A-B)_+}(\Delta)+ E_{A_+}(\Delta)+ E_{B_+}(\Delta).$$
	By (\ref{pogoj111}) this is equivalent to	$E_{(A+B)_++A_-+B_-}(\Delta) = E_{(A+B)_-+A_++B_+}(\Delta)$, which is
	true due to $(A+B)_++A_-+B_-=(A+B)_-+A_++B_+$.
	
	\textsl{Homogeneity of $m^\Delta$.} To prove $m^\Delta(\lambda A)=\lambda m^\Delta(A)$ for $A\in (W_1)_h, \lambda\in \RR$ it suffices to consider 
	$A\in (W_1)_+$
	(due to $\lambda A= \lambda A_+ - \lambda A_-$ and additivity of $m^\Delta$). For $\lambda\geq 0$ this is (\ref{pogoj222}).
	For $\lambda< 0$ we have
		$$m^\Delta(\lambda A) := E_{(\lambda A)_+}(\Delta)-E_{(\lambda A)_-}(\Delta)= 
		-E_{\left|\lambda\right| A}(\Delta)\underbrace{=}_{\text{by}\;(\ref{pogoj222})}
		-\left|\lambda\right| E_{A}(\Delta)=\lambda m^\Delta(A).$$
			
	\textsl{Boundedness of $m^\Delta$.} For $A\in (W_1)_h$ we have
				\beqn
					\left\|m^\Delta(A)\right\| &   =  & \left\|E_{A_+}(\Delta)-E_{A_-}(\Delta)\right\|
															   \leq  \left\|E_{A_+}(\Delta)\right\|+ \left\|E_{A_-}(\Delta)\right\|\\
															  & \underbrace{\leq}_{\text{by}\;(\ref{pogoj333})} & k_\Delta \left(\left\|A_+\right\| + 							
															  	\left\|A_-\right\|\right)
															  	\underbrace{\leq}_{\left\|A_+\right\|, \left\|A_-\right\|\leq \left\|A\right\|} 
															  2\left\|A\right\|k_\Delta.
			  \eeqn
\eprf


\section{Extension of the integration to $\overline{\cI(m)\otimes_\FF W_1}$}

Assume the notation from Section \ref{razdelek-z-nenegativnimi-merami}. 
Let $\left(X,\cS,W_1\subseteq B(\cH),W_2\subseteq B(\cK),m\right)$ be a measure space with a measure $m$, $\cI(m)$ the set of $m$-integrable functions and
$\FF\in\left\{\RR,\CC\right\}$. We equip
$\cI(m)\otimes_\FF W_1$ with a supremum norm, i.e., for every $F\in\cI(m)\otimes_\FF W_1$ we define 
$\left\|F\right\|_\infty:=\sup_{x\in X}\left\|F(x)\right\|$.
Let $\overline{\cI(m)\otimes_\FF W_1}$ be a norm completion of $\cI(m)\otimes_\FF W_1$.
For every $F\in \overline{\cI(m)\otimes_\FF W_1}$ we define
	\bequ \label{razsirimo-integracijo} \int_X F\; dm := \lim_{i\to\infty} \int_X F_i\;dm, \eequ
where $\{F_i\}_i$ is any sequence of elements from $\cI(m)\otimes_\FF W_1$ converging to $F$ in the supremum norm. 

The definition is well-defined by the following proposition.

\bpro \label{razsiritev-integracije-trditev}
	The integral $\int_X F\; dm$ exists and is independent of the choice of the sequence $\{F_i\}_i$.
\epro

\bprf
	Since $\{F_i\}_i$ is a Cauchy sequence, for each $\epsilon>0$ there exists $n_{\epsilon} \in\NN$ such that 
	$\left\|F_m-F_n\right\|_\infty<\epsilon$ for every $m, n\geq n_{\epsilon}$.  
	By (\ref{neenakost1}), 
	$\int_X \epsilon \cdot \id_{\cH}\;dm\succeq$ $\int_X (F_n-F_m)\;dm \succeq$ $-\int_X \epsilon \cdot \id_{\cH}\;dm,$
	where $\id_\cH$ denotes the identity operator on $\cH$. 
	Hence
		$\left\|\int_X (F_n-F_m)\;dm\right\|\leq \left\|\int_X (\epsilon \cdot \id_{\cH})\;dm\right\| 
		=\epsilon\left\|m_{\id_{\cH}}(X)\right\|.$
	Therefore  $\int_X F_i\; dm$ is a Cauchy sequence and hence convergent. 
	
	Independence of $\int_X F\; dm$ of the sequence is proved similarly. Namely, for
	the sequences $F_i$, $G_i$ converging to $F$, the sequence $F_i-G_i$ converges to $0$ and by the above
	argument $\int_X (F_i-G_i)\; dm$ is a convergent sequence with the limit 0. 
\eprf

For a locally compact Hausdorff space $X$ and a Banach space $\cY$ over $\FF\in\left\{\RR,\CC\right\}$,
let $C_c(X,\cY), C_0(X,\cY)$ be the vector spaces of $\cY$-valued functions with a compact support and $\cY$-valued functions
which vanish at infinity respectively, i.e., $F\in C_c(X,\cY)$ iff $F\in C(X,\cY)$ and the set $\left\{x\in X\colon F(x)\neq 0\right\}$ is compact
and $F\in C_0(X,\cY)$ iff $F\in C(X,\cY)$ and for every $\epsilon>0$ there exists a compact set $K_\epsilon$, such that $\left\|F(x)\right\|< \epsilon$ for every $x\in K_\epsilon^c$. 

Let $(X,\B(X),W_1,W_2,m)$ be a space with a measure $m$.
Since $C_c(X,\FF)\otimes_\FF W_1$ is dense in $C_0(X,W_1)$ endowed with the supremum norm (see \cite[Proposition 44.2.]{Treves}), we have
	$$C_0(X,W_1)\subseteq \overline{\cI(M)\otimes_\FF W_1}.$$
	
Let $V\leq C_0(X,W_1)$ be a vector subspace of $C_0(X,W_1)$ and $L:V\to W_2$ a bounded linear map. We call $L$ \textsl{positive} if 
$L(V_+)\subseteq (W_2)_+$, where $V_+:=V\cap C_0(X,W_1)_+$ is a positive cone of $V$ inherited from the positive cone
	$$C_0(X,W_1)_+:=\left\{F\in C_0(X,W_1) \colon F(x)\in (W_1)_+ \text{ for every } x\in X\right\}$$ of $C_0(X,W_1).$

Theorem \ref{tretji-Rieszov-operatorski-izrek} is a version of the Riesz representation theorem and slightly extends \cite[Proposition 2]{Cim-Zal} 
from the case of a positive bounded linear map $L:C_c(X,\RR)\otimes_\RR B(\cH)_h\to B(\cK)_h$ on a locally compact and $\sigma$-compact metrizable space $X$ and real Hilbert spaces $\cH, \cK$, to the case of a positive bounded linear map $L: C_0(X,W_1) \to W_2$ on a locally compact Hausdorff space $X$ and Hilbert spaces $\cH$,
$\cK$ over $\FF\in\left\{\RR,\CC\right\}$.

\bthe	\label{tretji-Rieszov-operatorski-izrek}
	Let $X$ be a locally compact Hausdorff space, $\cH$, $\cK$ Hilbert spaces over $\FF\in \left\{\RR,\CC\right\}$ and 
	$W_1\subseteq B(\cH)$, $W_2\subseteq B(\cK)$ von Neumann algebras.
	\begin{enumerate}
		\item $\FF=\RR$: For every positive bounded linear map $L: C_0(X,(W_1)_h) \to (W_2)_h$ there exists a unique regular non-negative measure
				$$m: \Bor(X)\rightarrow B((W_1)_h, (W_2)_h)$$ 
			such that $L(F)=\int_X {F\; dm}$ holds for all $F\in C_0(X,(W_1)_h).$
		\item $\FF=\CC$: For every positive bounded linear map $L: C_0(X,W_1) \to W_2$ there exists a unique regular non-negative measure
				$$m: \Bor(X)\rightarrow B(W_1, W_2)$$ 
			such that $L(F)=\int_X {F\; dm}$ holds for all $F\in C_0(X,W_1).$
	\end{enumerate}
\ethe

\bprf
	We replace some assumptions of \cite[Proposition 2]{Cim-Zal} by weaker ones stepwise:
	\begin{enumerate}
	\item \textsl{Replacing a locally compact and $\sigma$-compact metrizable space $X$
	  by a locally compact Hausdorff space $X$:}
	  \begin{enumerate}
	  	\item $\FF=\RR$: Since $C_c(X,\RR)$ is dense in $C_0(X,\RR)$, $T$ in \cite[Theorem 19]{Ber} can be uniquelly extended to the bounded map on
				$C_0(X,\RR)$. By complexifying $\cH$, $\cK$ and linearly extending $T$ to the complexification $C_0(X,\CC)$ of $C_0(X,\RR)$,
				$T$ remains a positive bounded linear map.
				The construction of the representing measure $E$ of $T$ in the proof of \cite[Theorem 19]{Ber} remains the same, just that we use the
				version of Riesz theorem for $C_0(X,\CC)$ (see \cite[C.17.~Theorem.]{Con}), use Lemma \ref{density} and work with nets
				as in the proof of \cite[IX.1.14.~Theorem.]{Con}. Applying this to \cite[Proposition 2]{Cim-Zal} yields (1).
	  	\item $\FF=\CC$: Restricting a positive bounded linear map $L:C_c(X,\CC)\otimes_\CC B(\cH)\to B(\cK)$ to 
	  		a positive bounded linear map $L:C_c(X,\RR)\otimes_\RR B(\cH)_h\to B(\cK)_h$ and applying the proof of $\FF=\RR$ case above yields the statement of
	  		\cite[Proposition 2]{Cim-Zal} for the locally compact Hausdorff space and $L:C_c(X,\CC)\otimes_\CC B(\cH)\to B(\cK)$.
	  \end{enumerate}
	\item \textsl{Replacing $B(\cH), B(\cK)$ by $W_1,W_2$:} This follows trivially.	
	\item \textsl{Replacing $C_c(X,\FF)\otimes_\FF W_1$ by $C_0(X,W_1)$:}
	By (2), $L|_{C_c(X,\FF)\otimes_\FF W_1}$ has a unique non-negative representing measure 
	$m$. For $F\in C_0(X,W_1)$ there exists a sequence 
	$F_i\in C_c(X,\FF)\otimes_\FF W_1$ such that $\lim_i F_i=F$ and by the continuity of $L$ we have $L(F)=\lim_i L(F_i)$.
	By the definition also $\int_X F\;dm=\lim_{i} \int_X F_i\;dm$.
	Hence $L(F)=\int_X F\;dm$.
	\end{enumerate}
\eprf



\section{Connection with measures and integration from \cite{Mitt-Young}}

Let $(X,\cS,W_1\subseteq B(\cH),W_2\subseteq B(\cK),m)$ be a space with a measure $m$ (see Section \ref{razdelek-z-nenegativnimi-merami}).
From now on up to the end of this section we will assume that $X$ is a locally compact Hausdorff space, $\cS=\B(X)$ is a $\sigma$-algebra of Borel sets on $X$
and $m_A$ are regular measures for every $A\in (W_1)_+$, i.e., for every $A\in (W_1)_+$ and every $k_1,k_2\in\cK$
	$$\left\langle m_A(\cdot)k_1,k_2\right\rangle:\B(X)\to \CC,\quad \Delta\mapsto \left\langle m_A(\Delta)k_1,k_2\right\rangle$$ 
is a regular complex measure. 

\textsl{Semivariation} of $m$ is the map $\overline{m}:\Bor(X)\to [0,\infty]$ defined by
	$$\overline{m}(X) := \sup\left\{\left\|\sum_{j=1}^n m_{A_j}(\Delta_j)\right\|\right\},$$
where the supremum is taken over all finite collections of disjoint sets $\Delta_1,\Delta_2,\ldots,\Delta_n$ such that $X=\cup_{j=1}^n \Delta_j$
and all $A_1,A_2,\ldots,A_n\in W_1$ with norm at most $1$.

Let $(W_2)_\ast$ denote a predual of a von Neumann algebra $W_2$. Recall that $(W_2)_\ast$ is the set of all ultra-weakly (or equivalently ultra-strongly) continuous linear functionals on $W_2$ (see \cite[I.3.~Theorem 1.(iii)]{Dix}). 
For every $T\in (W_2)_\ast$ and every $A\in W_1$ we define a map 
	$$\left\langle T,m_A(\cdot)\right\rangle:\B(X)\to \CC,\quad \left\langle T,m_A(\Delta)\right\rangle:= T(m_A(\Delta)).$$ 
		
The next proposition shows that by the assumptions of the first paragraph our non-negative measures have the properties of measures obtained in \cite[Theorem 3.3]{Mitt-Young}.

\bpro \label{opomba-Mitter-Young}
			Assume $X$, $\B(X)$ and $m_A$ are as above.
			We claim that $m$ is a finitely additive measure with a finite semivariation, i.e., $\overline{m}(X)<\infty$, such that for every $T\in (W_2)_{\ast}$
			and every $A\in W_{1}$ the maps $\left\langle T,m_A(\cdot)\right\rangle$ are regular countably additive complex measures with a bounded variation.
\epro
			\begin{proof}
				Firstly we prove $\overline{m}(X)<\infty.$
				\beqn
					&&\overline{m}(X) = \sup\left\{\left\|\sum_{j=1}^n m_{A_j}(\Delta_j)\right\|\right\}\\
					&=& \sup\left\{\left\|\sum_{j=1}^n (m_{\re(A_j)}(\Delta_j) - i\cdot m_{\im(A_j)}(\Delta_j))\right\|\right\}\\
					&=& \sup\left\{\left\|\sum_{j=1}^n ((m_{\re(A_j)_+}- m_{\re(A_j)_-})(\Delta_j) + 
						i (m_{\im(A_j)_+} - m_{\im(A_j)_-}(\Delta_j)))\right\|\right\}\\
					&\leq&  4\sup\left\{\left\|\sum_{j=1}^n m_{\id_\cH}(\Delta_j)\right\|\right\}	= 4\left\|m_{\id_\cH}(X)\right\|<\infty,
				\eeqn
				where the supremum is taken over all finite collections of disjoint sets $\Delta_1,\Delta_2,\ldots$, $\Delta_n$, such that $X=\cup_{j=1}^n \Delta_j$
				and all $A_1,A_2,\ldots,A_n\in W_1$ with norm at most $1$, and $\id_\cH$ denotes the identity operator on a Hilbert space $\cH$.
				Note that the first inequality in the last line follows by $\left\|\re(A_j)_\pm\right\|$, $\left\|\im(A_j)_\pm\right\|\leq 1$, which further implies
				$\re(A_j)_\pm$, $\im(A_j)_\pm \preceq \id_{\cH}$
				and finally $m_{\re(A_j)_\pm}(\Delta_j)$, $m_{\im(A_j)_\pm}(\Delta_j) \preceq m_{\id_\cH}(\Delta_j)$ for every $j=1,2,\ldots,n$.
				
				We have to prove that the measures 
					\begin{equation}\label{equality-9}
						\left\langle T,m_A(\cdot)\right\rangle: \B(X)\to \CC, \quad \Delta \mapsto \left\langle T,m_A(\Delta)\right\rangle,
					\end{equation}
				where $T\in (W_2)_\ast$ and $A\in W_1$, are countably additive, regular and have a finite variation. By the usual decompositions of $T$ and $A$ into 
				the linear combination of the positive elements it suffices to take $T\in ((W_2)_\ast)_+$ and $A\in (W_1)_+$ (for $T$ this follows by \cite[I.4.~Theorem 6.(i)]{Dix}). 
				Here $((W_2)_\ast)_+$ denotes the set of all $T\in (W_2)_\ast$, such
				that $T(B)\geq 0$ for every $B\in (W_2)_+$.
				
				Let us first prove the countable additivity.
				Take $\Delta\in \B(X)$ and a sequence $\Delta_j \subseteq \B(X)$, such that $\Delta=\cup_{j=1}^\infty \Delta_j$ and $\Delta_j$ are mutually disjoint sets.
				To prove 
						$\left\langle T,m_A(\Delta)\right\rangle=\sum_{j=1}^{\infty} \left\langle T,m_A(\Delta_j)\right\rangle$
				it suffices to show that
				$m_A(\Delta)=\sum_{j=1}^{\infty} m_A(\Delta_j)$ in the ultra-strong topology of $W_2$. 
				Since the ultra-strong and strong topologies coincide on bounded subsets of $W_2$, this imediatelly follows by the definition of $m_A$.
				
				By the countable additivity, it follows that the sum $\sum_{j=1}^{\infty} \left\langle T,m_A(\Delta_j)\right\rangle$ is absolutely convergent and hence
				the variation of (\ref{equality-9}) is finite.
				
				To prove the inner regularity, let us take an open set $U$. 
				Since $\left\langle T,m_A(\cdot)\right\rangle$ is a finite positive
				measure, it suffices to prove that for every $\epsilon>0$ there exists a compact set $K_\epsilon\subset U$, such that 
				$\left\langle T,m_A(U\setminus K_\epsilon)\right\rangle<\epsilon$. 
				Since $T$ is ultra-strongly continuous,
				there exists a sequence $\{k_j\}_{j=1}^\infty\subseteq \cK$ 
				with 		
				$\sum_{j=1}^{\infty}\left\|k_j\right\|^2<\infty$, such that for every $B\in W_2$ satisfying $\sum_{j=1}^\infty\|(m_A(U)-B)k_j\|^2<1$,
				it follows $\left|\left\langle T,m_A(U)-B\right\rangle\right|<\epsilon$.
				There exists $N\in\NN$, such that $\sum_{j=N+1}^{\infty}\left\|k_j\right\|^2<\frac{1}{2^2\left\|m_A(X)\right\|^2}$.
				By the inner regularity of the measures $m_A$, for every $\ell=1,2,\ldots,N$ there exists a compact set $K_\ell\subset U$, such that
					\begin{equation}\label{enacba}
						\left\langle m_A(U\setminus K_\ell)k_\ell,k_\ell\right\rangle< \frac{1}{2^2 N^2\left\|m_A(X)\right\|^3 \left\|k_\ell\right\|^2}.
					\end{equation}
				Therefore
					\beqn
						{\left\|m_A(U\setminus K_{\ell})k_\ell\right\|}^4&=& 
						 {\left\langle m_A(U\setminus K_{\ell})k_\ell, m_A(U\setminus K_\ell)k_\ell\right\rangle}^2 \\ 
						&\leq& \left\langle m_A(U\setminus K_{\ell})k_{\ell}, k_{\ell}\right\rangle\cdot \left\langle \left(m_A(U\setminus K_{\ell})\right)^2k_{\ell},  m_A(U\setminus 	
							K_{\ell})k_{\ell}\right\rangle \\
							&\leq&  \left\langle m_A(U\setminus K_{\ell})k_{\ell}, k_{\ell}\right\rangle \cdot {\left\|m_A(U\setminus K_{\ell})\right\|}^3 {\left\|k_{\ell}\right\|}^2\\
							&\leq&  \left\langle m_A(U\setminus K_{\ell})k_{\ell}, k_\ell\right\rangle \cdot {\left\|m_A(X)\right\|}^3 {\left\|k_\ell\right\|}^2\\
							&<&  \frac{1}{2^2N^2},
					\eeqn
				where we used Cauchy-Schwarz inequality for the semi-inner product 
				$[k_1,k_2]:=\left\langle m_A(U\setminus K_\ell)k_1,k_2\right\rangle,$ $k_1,k_2 \in \cK$ in the first inequality (Notice that $m_A$ is a positive operator-valued measure, since we are considering $A$ from $(W_1)_+$.), 
				Cauchy-Schwarz inequality in $\cK$ in the second, finite additivity of $m_A$ in the third (i.e., from $U\setminus K_{\ell}\subseteq X$ it follows 
				$0\preceq m_A(U\setminus K_{\ell})\preceq m_A(X)$ and hence $\left\|m_A(U\setminus K_{\ell})\right\|\leq \left\|m_A(X)\right\|$), and (\ref{enacba}) in the last one.
				Hence, for $K_\epsilon:=\cup_{\ell=1}^N K_\ell$ it follows $\sum_{\ell=1}^N \left\|m_A(U\setminus K_\epsilon)k_\ell\right\|^2\leq \frac{1}{2}.$
				Therefore
					\beqn 
									 \sum_{\ell=1}^\infty \left\|m_A(U\setminus K_\epsilon)k_\ell\right\|^2
						&<& \frac{1}{2} + \left\|m_A(U\setminus K_\epsilon)\right\|^2 \sum_{\ell=N+1}^{\infty} \left\|k_\ell\right\|^2\\
						&\leq&   \frac{1}{2} + \left\|m_A(X)\right\|^2 \frac{1}{2\left\|m_A(X)\right\|^2}
						\leq 	 1.
					\eeqn
				It follows that $\left\langle T,m_A(U\setminus K_\epsilon)\right\rangle<\epsilon$.
				Since $\epsilon>0$ was arbitrary, this proves the inner regularity.
				
				The outer regularity is proved analoguously.
				\end{proof}
By Proposition \ref{opomba-Mitter-Young}, $m$ has a finite semivariation. Now we will compare our integration with the integration with respect
to the measure of finite semivariation from \cite[Section 3]{Mitt-Young}. Let $\cB$ stand for the simple Borel measurable functions on $X$.
By (\ref{enakost3}), Theorem \ref{konvergencni-izrek} and (\ref{razsirimo-integracijo}), it easily follows that the integrations coincide on $\overline{\cB\otimes W_1}$. But $\overline{\cB\otimes W_1}$ consists just of bounded $W_1$-valued function, while not all elements from $\overline{\cI(m)\otimes W_1}$ are necessarily bounded. Hence, our integration theory extends the integration theory, when $m$ is regarded as a finitely additive measure with a finite semivariation.

Now we will comment on the connection between Theorem \ref{tretji-Rieszov-operatorski-izrek} and \cite[Theorem 3.3.]{Mitt-Young}.
Since every von Neumann algebra is a dual of a Banach space, $L: C_0(X,W_1) \to W_2$ from Theorem \ref{tretji-Rieszov-operatorski-izrek} 
has a representing measure given by \cite[Theorem 3.3.]{Mitt-Young}. By Proposition \ref{opomba-Mitter-Young} and the uniqueness of the measures from 
Theorem \ref{tretji-Rieszov-operatorski-izrek} and \cite[Theorem 3.3.]{Mitt-Young}, it follows that both representing measures coincide. 
Hence, we derived another proof for the special case of \cite[Theorem 3.3.]{Mitt-Young} and obtained a concrete description of measures given by \cite[Theorem 3.3.]{Mitt-Young} for this special case.


\section{Non-negative spectral measures}
\label{nenegativne-spektralne}

Notation remains as in Section \ref{razdelek-z-nenegativnimi-merami}. 
Let $(X,\cS,W_1\subseteq B(\cH),W_2\subseteq B(\cK))$ be a measure space, where $\cH$, $\cK$ are complex Hilbert spaces.
We denote by $(W_1)_p$ the set of all hermitian projections in $W_1$, i.e., $A\in (W_1)_p \Leftrightarrow A=A^\ast=A^2$.
	 
	A non-negative measure  
		$$M \colon \cS \to B(W_1,W_2)$$ 
	is \textit{a non-negative spectral measure} if for every $P\in (W_1)_p$ the set function
		$$M_P \colon \cS \to W_2, \quad M_P(\Delta):=M(\Delta)(P),$$
	is a spectral measure and if the equality 
		$$M_P(\Delta_1)M_Q(\Delta_2)=M_{PQ}(\Delta_1\cap \Delta_2) $$
	holds for all hermitian projections $P,Q\in (W_1)_p$ and all sets $\Delta_1,\Delta_2\in \cS$.
A pentuple $(X,\cS,W_1,W_2,M)$ is \textsl{a space with a measure $M$}.

\begin{remark}
	\begin{enumerate}
		\item Recall that in the introduction we defined $M$ to be just \emph{a set function} with the above properties and not \emph{a non-negative measure}. However,
			in the next section (see Corollary \ref{prvotna-definicija}) we will show that a set function with the above properties is automatically
			a non-negative measure.
		\item We will need non-negative spectral measures to represent $\ast$-representations of the form $\rho: C(X,W_1) \to W_2$ (see Section 
			\ref{integralske-reprezentacije}). Since in the case $\cH$ is a real Hilbert space, $C(X,W_1)$ is not even an algebra, we cannot study $\ast$-representations
			$\rho$. Therefore, from now on all the Hilbert spaces will be complex.
	\end{enumerate}
	
\end{remark}

By the well-known result, e.g., \cite[Theorem 5.1]{Sch2012}, every hermitian operator $A\in (W_1)_h\subseteq B(\cH)_h$, has a unique spectral measure $E:\B([a,b])\to B(\cH)$ such that 
	$A=\int_{[a,b]} \lambda\;dE(\lambda)$ and $\sigma(A)\subseteq [a,b]$, where $\sigma(A)$ denotes the spectrum of $A$.
From $A\in (W_1)_h$ it follows by \cite[I.~Proposition 2]{Dix}, that $E(\Delta)\in W_1$ for every $\Delta\in \B([a,b])$.
By \cite[p.~63-64]{Sch2012}, there exists a sequence $S_\ell(A)$ of Riemann sums of the form
	\bequ \label{Riemannove-vsote-operatorja} S_\ell(A)=\sum_{k=1}^{n_{\ell}} \zeta_{k,\ell} 	
		(E(\lambda_{k,\ell})-E(\lambda_{k-1,\ell})), \eequ
where the family $\left\{E(\lambda)\mid \lambda\in\RR\right\}$ is the resolution of identity (see \cite[Definition 4.1]{Sch2012}) corresponding to the spectral measure $E$, 
$S_{\ell+1}(A)$ is a refinement of $S_\ell(A)$ 
and $\left\|A-S_\ell(A)\right\|\leq \frac{1}{\ell}.$

Recall that the set $\cI(M)$ of all \textsl{$M$-integrable} functions is a complex vector space and 
it consists of at least all bounded $\cS$-measurable complex functions (Here we integrate with respect to $M$ by regarding $M$ as a non-negative measure and use integration theory from Section \ref{razdelek-z-nenegativnimi-merami}.).

The next proposition shows that the integration with respect to a non-negative spectral measure is multiplicative.

\bpro \label{integracija-po-multiplikativni-meri-je-mutliplikativna}
	Let $(X,\cS,W_1,W_2,M)$ be a space with a non-negative spectral measure $M$ and $F,G$ elements from $\overline{\cI(M)\otimes W_1}$. Then the equality
		\bequ \label{multiplikativnostna-enakost} \int_X FG\; dM=\left(\int_X F\; dM\right)\left(\int_X G\; dM\right) \eequ
	holds.
\epro

\bprf
	We will prove (\ref{multiplikativnostna-enakost}) in two steps. First we will consider the case $F,G\in \cI(M)\otimes W_1$ and then use it in the proof of the 
	general case $F,G\in \overline{\cI(M)\otimes W_1}$.
	
	\textsl{Case 1 - $F,G\in \cI(M)\otimes W_1$:}
	By the linearity it suffices to consider $F=f\otimes A, G=g\otimes B$ for $f,g\in \cI(M)_+$, $A,B\in (W_1)_+$ and 	
	prove
		\bequ \label{multiplikativnost-na-projektorjih}
			\int_X \left(fg \otimes AB\right) \; dM =
			\left(\int_X (f\otimes A) \; dM\right)\left(\int_X (g \otimes B) \; dM\right)
		\eequ
		
	Let us first show that (\ref{multiplikativnost-na-projektorjih}) holds 
	for every $A,B\in (W_1)_p$. 
	There are increasing sequences $\{s_k\}_{k}$, $\{t_k\}_{k}$ of simple functions such that
		$\lim_{k} s_k=f, \lim_{k} t_k=g$. 
	By Theorem \ref{konvergencni-izrek}, 
		\beqn
			\lim_{k}\int_X \left(s_k t_k \otimes AB\right) \; dM &=& \int_X \left(fg \otimes AB\right) \; dM, \\
			\lim_{k}\int_X \left(s_k \otimes A\right) \; dM &=& \int_X \left(f \otimes A\right) \; dM, \\
			\lim_{k}\int_X \left(t_k \otimes B\right) \; dM &=& \int_X \left(g \otimes B\right) \; dM,
		\eeqn
	where all the limits are in the strong operator topology (For the first equality we have also
	used the decomposition of $PQ$ into four positive parts and applied the convergence theorem to each of them.).
	By the definition of $M$ and the linearity of the integration, the equality (\ref{multiplikativnost-na-projektorjih}) is true for all simple functions $s_k,t_k$ and for every $A,B\in (W_1)_p$. Hence 
	it is true also for every $f,g\in \cI(M)_+$ and every $A,B\in (W_1)_p$.
	
	Let now $A,B\in (W_1)_+$ be arbitrary. If $S_\ell(A)$, $S_\ell(B)$ are defined as in (\ref{Riemannove-vsote-operatorja}), then  
		$$\frac{1}{\ell}\id_\cH \succeq A-S_\ell(A)\succeq 0, \quad
		\frac{1}{\ell}\id_\cH \succeq B-S_\ell(B)\succeq 0,$$ 
	where $\id_\cH$ denotes the identity operator of $\cH$, 
	and by (\ref{neenakost1}) we get
		$$\left\|\int_X (f\otimes (A-S_\ell(A)))\; dM\right\| \leq \frac{1}{\ell} \left\|\int_X (f\otimes \id_\cH)\; dM\right\|,$$
		$$\left\|\int_X (g\otimes (B-S_\ell(B)))\; dM\right\| \leq \frac{1}{\ell} \left\|\int_X (g\otimes \id_\cH)\; dM\right\|.$$
	Therefore 
		\beqn
			\int_X (f\otimes A)\; dM &=& \lim_{\ell}\int_X (f\otimes S_\ell(A))\; dM\\
			\int_X (g\otimes B)\; dM &=& \lim_{\ell}\int_X (g\otimes S_\ell(B))\; dM.
		\eeqn
	Since (\ref{multiplikativnost-na-projektorjih}) holds
	for all hermitian projections $A, B\in (W_1)_p$, it follows by the linearity that
	\beqn
			\int_X (fg\otimes S_\ell(A)S_\ell(B))\; dM &=& 
			\left(\int_X (f\otimes S_\ell(A))\; dM \right)\left( \int_X (g\otimes S_\ell(B))\; dM\right).
	\eeqn
	To prove (\ref{multiplikativnost-na-projektorjih}) for $A,B\in (W_1)_+$, 
	it remains to prove that
		\bequ 
			\label{enakost-100}\int_X (fg\otimes AB)\; dM = \lim_{\ell\to \infty} \int_X (fg\otimes S_\ell(A)S_\ell(B))\; dM. 
		\eequ
	We denote $C_\ell:=AB-S_\ell(A)S_\ell(B)$ and $\epsilon_\ell:=\left\|C_\ell\right\|$.
	By the usual decomposition of $C_\ell$ into the linear combination of four positive elements, we conclude that
			$$\left\|\int_X \left(fg\otimes C_\ell\right)\; dM\right\| 	
			\leq 4\epsilon_\ell \left\| \int_X \left(fg\otimes \id_\cH\right)\; dM\right\|.$$	
	Here we used $\left\|\re(C_\ell)_{\pm}\right\|, \left\|\im(C_{\ell})_{\pm}\right\|\leq \left\|C_\ell\right\|$
	and (\ref{neenakost1}). Since $\lim_{\ell}\epsilon_\ell=0$, (\ref{enakost-100}) follows.
	
	\textsl{Case 2 - $F,G\in \overline{\cI(M)\otimes W_1}$:}
	By the definition of the integration with respect to $M$ (see (\ref{razsirimo-integracijo})),
		$$\int_X F\; dM=\lim_{i}\int_X F_i\; dM, \quad
		\int_X G\; dM=\lim_{i}\int_X G_i\; dM$$
	where $F_i, G_i\in \cI(M)\otimes W_1$ are any sequences converging to $F$, $G$ in the 
	supremum norm. 
	Since $F_iG_i$ converges to $FG$ in the supremum norm, it follows
		$\int_X FG\; dM=\lim_{i}\int_X F_iG_i\; dM$.
	By Case 1, equality (\ref{multiplikativnostna-enakost}) holds for every pair $F_i, G_i$; hence also for the pair
	$F,G$. 
\eprf

\brem
	In the sequel we will use the statement of Proposition \ref{integracija-po-multiplikativni-meri-je-mutliplikativna} just for the
	pairs $F, G$ from the set $C_0(X,W_1)$, where $X$ is a locally compact Hausdorff space $X$ and $\B(X)$ a Borel $\sigma$-algebra on $X$. 
	Since elements from $C_0(X,W_1)$ are bounded, it would suffice to prove the validity of the statement of Proposition \ref{integracija-po-multiplikativni-meri-je-mutliplikativna} 
	in a much lesser generality, i.e., bounded elements $F,G$ from $\overline{\cI(M)\otimes W_1}$ would do the job.
	But this is simple. Let us write it down:
	
	Take bounded elements $F,G$ from $\overline{\cI(M)\otimes W_1}$.
	There are sequences $\{S_k\}_k$, $\{T_k\}_k$ of simple functions from $\cI(M)\otimes W_1$ (i.e., $S_k=\sum_i s_{ki}\otimes A_{ki}$ and 
	$T_k=\sum_j t_{kj}\otimes B_{kj}$, where $s_{ki}, t_{kj}$ are the usual simple functions and $A_{ki}, B_{kj}\in W_1$),
	such that 
		$\lim_{k} \left\|F-S_k\right\|_\infty=0$ and $\lim_k \left\|G-T_k\right\|_\infty=0$.
	Therefore it is also true that
		$\lim_{k} \left\|FG-S_kT_k\right\|_\infty=0$.
	By the definition of the integration with respect to $M$ (see (\ref{razsirimo-integracijo})), 
		$\int_X F\;dM$$=\lim_k \int_X S_k\;dM$, 
		$\int_X G\;dM$$=\lim_k \int_X T_k\;dM$ and
		$\int_X FG\;dM$$=\lim_k \int_X S_kT_k\;dM$.
	By the linearity and the multiplicativity of $M$ on hermitian projections, it folows that 
	$\left(\int_X S_k\;dM\right)\left(\int_X T_k\;dM\right)=$$\int_X S_kT_k\;dM$.
	Therefore (\ref{multiplikativnostna-enakost}) holds for $F, G$.
\erem

\section{Characterization of non-negative spectral measures}
\label{razdelek-karakterizacija-spektralnih-mer}

Let $(X,\cS,W_1,W_2)$ be a measure space (see Section \ref{razdelek-z-nenegativnimi-merami}). In Theorem \ref{karakterizacija-nenegativnih-kompleksnih-mer} we 
characterized non-negative measures on $(X,\cS,W_1,W_2)$ via families $\left\{E_A\right\}_{A\in (W_1)_+}$ of positive operator-valued measures. We would like to have an analoguous characterization in the case of non-negative spectral measures.

The following theorem provides a characterization of non-negative spectral measures on $(X,\cS,W_1,W_2)$ via families $\left\{F_P\right\}_{P\in (W_1)_p}$ of spectral measures (Recall that $(W_1)_p$ denotes the set of all hermitian projections in $W_1$.).
This characterization will be used to prove our main results (see Theorem \ref{reprezentacijski-izrek-1} and Corollary \ref{reprezentacijski-izrek-2}) in the next section. 

\bthe\label{karakterizacija-nenegativnih-spektralnih-mer}
	Let $(X,\cS,W_1,W_2)$ be a measure space, $\left\{F_P\right\}_{P\in (W_1)_p}$ a family of 
	spectral measures $F_P:\cS\to W_2$.
	
	There is a unique non-negative spectral measure $M$ such that
		$$M_P=F_P$$
	for all hermitian projections $P\in (W_1)_p$ iff the following conditions hold.
		\bequ\label{pogoj1111}  \sum_{i=1}^{n}\lambda_i F_{P_i}(\Delta)=\sum_{j=1}^m \mu_j F_{Q_j}(\Delta),\eequ
	for all hermitian projections $P_i,Q_j\in (W_1)_p$, all real numbers $\lambda_i,\mu_j \in\RR$, and all sets $\Delta\in\cS$ such that
		$\sum_{i=1}^{n}\lambda_i P_i=\sum_{j=1}^m \mu_j Q_j$,
	for each set $\Delta\in\cS$ there exists a constant $k_\Delta\in\RR^{>0}$ such that  
		\bequ\label{pogoj3333} \left\|F_P(\Delta)\right\|\leq k_\Delta \eequ
	for all hermitian projections $P\in (W_1)_p$, and
		\bequ\label{pogoj4444} M_P(\Delta_1)M_Q(\Delta_2)=M_{PQ}(\Delta_1\cap \Delta_2) \eequ
	holds for all hermitian projections $P,Q\in (W_1)_p$ and all sets $\Delta_1,\Delta_2\in \cS.$
\ethe

Every family $\left\{F_P\right\}_{P\in (W_1)_p}$ which satisfies the conditions above is called 
\textsl{a compatible family of spectral measures}.

\bprf
	The nontrivial direction is the if part. Suppose that we are given a family $\left\{F_P\right\}_{P\in (W_1)_p}$ of 
	spectral measures $F_P:\cS\to W_2$, which satistfies the conditions (\ref{pogoj1111}), (\ref{pogoj3333}), (\ref{pogoj4444}). By the statement of Theorem 
	\ref{karakterizacija-nenegativnih-spektralnih-mer}, we have to find a non-negative measure $M:\cS\to B(W_1,W_2)$ such that $M_P=F_P$ for all $P\in (W_1)_p$, 
	where $M_P:\cS\to W_2$ is defined by $M_P(\Delta):=M(\Delta)(P)$ for every $\Delta\in\cS$. 
	Therefore all that remains is to define the set functions $M_A\colon \cS \to W_2$
	for every $A\in (W_1)_+$ and prove that the family $\{M_A\}_{A\in (W_{1})_+}$ is a well-defined family of positive operator-valued measures, which 
	satisfies the conditions (\ref{pogoj111}), (\ref{pogoj222}), (\ref{pogoj333}) of Theorem \ref{karakterizacija-nenegativnih-kompleksnih-mer}. 
	Take $A\in (W_1)_+$. We separate two possibilities:
		\benu
			\item[(i)] If $A$ has a finite spectral decomposition $\sum_{k=1}^n\lambda_k P_k$, where
				$P_k$ are mutually orthogonal hermitian projections (i.e., $P_iP_j=0$ for every $i\neq j$), then 
					$$M_A(\Delta):= \sum_{k=1}^n\lambda_k F_{P_k}(\Delta).$$
			\item[(ii)] If $A$ does not have a finite spectral decomposition, then for $S_\ell(A)$ as in (\ref{Riemannove-vsote-operatorja})
					$$M_A(\Delta):= \lim_{\ell} M_{S_\ell(A)}(\Delta),$$ 
				where the limit is taken in the norm topology.
		\eenu
	We will prove the facts we need stepwise:
	
	\textsl{Step 1 - existence and uniqueness of $M_A(\Delta)$:} For $A$ with a finite spectral decomposition both facts are clear.
		The tougher part is to prove them for $A$ without a finite spectral decomposition. 
		Take the sequence $S_\ell(A)$ as in (ii). By the definition of $S_\ell(A)$, we conclude that 
		$\left\|S_{\ell_1}(A)-S_{\ell_2}(A)\right\|\leq \frac{2}{\ell}$ for $\ell_1,\ell_2\geq \ell$ and
		$S_{\ell_1}(A)-S_{\ell_2}(A)$ has a finite spectral decomposition. Let us denote it by $\sum_k\lambda_{k} Q_{k}$, where $Q_k$ are mutually orthogonal 
		hermitian projections and 
		$\lambda_k$ real numbers. Hence
			\beqn 
				&&\left\|M_{S_{\ell_1}(A)}(\Delta)- M_{S_{\ell_2}(A)} (\Delta)\right\|	\underbrace{=}_{
				\substack{\text{by}\;(\ref{pogoj1111})}}
				\left\|\sum_k \lambda_{k} M_{Q_{k}}(\Delta)\right\|\\
				&\underbrace{\leq}_{\substack{M_{Q_{k}}(\Delta)\succeq 0,\\\left|\lambda_{k}\right|\leq \frac{2}{\ell}}}&	
					\frac{2}{\ell} \left\|\sum_k M_{Q_{k}}(\Delta)\right\|
					\underbrace{=}_{\text{by}\;(\ref{pogoj1111})}
					\frac{2}{\ell}\left\|M_{\sum_k Q_{k}}(\Delta)\right\|
				\underbrace{\leq}_{\substack{\text{by}\;(\ref{pogoj3333})\\\text{for}\;\sum_k Q_{k}}} \frac{2}{\ell} k_\Delta
			\eeqn	
		(Note that for the last inequality $\sum_k Q_{k}$ has to be a hermitian projection, which is true since $Q_k$ are mutually orthogonal hermitian projections.).
		Therefore the sequence $M_{S_\ell(A)}(\Delta)$ is Cauchy in $W_2$ and hence convergent. So the operator $M_A(\Delta)$
		exists. Its uniqueness is proved analoguously, namely if $\tilde S_\ell(A)$ is another sequence satisfying (\ref{Riemannove-vsote-operatorja}),
		then $M_{S_{\ell}(A)}(\Delta)-M_{\tilde S_{\ell_1}(A)}(\Delta)$ is a Cauchy sequence converging to 0. 
		
		\textsl{Step 2 - $M_A$ is a positive operator-valued measure:} For $A\in (W_1)_+$ with a finite spectral decomposition this is clear. 
		Assume $A\in (W_1)_+$ does not have a finite spectral decomposition. Notice that all constants $\zeta_{k,\ell}$ in (\ref{Riemannove-vsote-operatorja}) can be 
		chosen such that $\zeta_{k,l}\geq 0$. Using this fact and by 			
			$$M_A(\Delta) \underbrace{=}_{(\ref{pogoj1111})} \lim_{\ell} \sum_{k=1}^{n_{\ell}} 
						\zeta_{k,\ell} (M_{E(\lambda_{k,\ell})-E(\lambda_{k-1,\ell})}(\Delta)), $$
		where $E(\lambda_{k,\ell})-E(\lambda_{k-1,\ell})\succeq 0$,
		it follows that $M_A(\Delta)\in (W_2)_+$.
		For $M_A$ to be a positive operator-valued measure we have to prove also the countable additivity. Take $\Delta=\cup_{j=1}^\infty \Delta_j$, where $\Delta,\Delta_j\in \cS$ and $\Delta_j$ are mutually disjoint. 
		Then the equality $M_A\left(\cup_{j=1}^\infty \Delta_j\right)= \sum_{j=1}^\infty M_{A}\left(\Delta_j\right)$ holds in the strong operator topology by tbe following:
		\beqn
			  M_{S_\ell(A)}\left(\cup_{j=1}^\infty \Delta_j\right)
			 &=&  \sum_{k=1}^{n_{\ell}} \zeta_{k,\ell} \left(M_{E(\lambda_{k,\ell})-E(\lambda_{k-1,\ell})}\left(\cup_{j=1}^\infty \Delta_j\right)\right) \\
			 &=& \sum_{i=1}^k\sum_{j=1}^\infty \zeta_{k,\ell} \left(M_{E(\lambda_{k,\ell})-E(\lambda_{k-1,\ell})}\left(\Delta_j\right)\right)
			 \\
			 &=&  \sum_{j=1}^\infty\sum_{k=1}^{n_{\ell}} \zeta_{k,\ell} \left(M_{E(\lambda_{k,\ell})-E(\lambda_{k-1,\ell})}\left(\Delta_j\right)\right)
			 = \sum_{j=1}^\infty M_{S_\ell(A)}(\Delta_j), 
		\eeqn		
	where the first and the forth equality hold by (\ref{pogoj1111}), the second by $M_{E(\lambda_{k,\ell})-E(\lambda_{k-1,\ell})}$ being spectral measures and 
	the third holds since all the operators in the sum are positive. Note that the second equality holds in the strong operator topology and not necessarily in the norm one, but this is all what we need in the proof.
	
	\textsl{Step 3 - $\{M_A\}_{A\in (W_{1})_+}$ satisfies the condition (\ref{pogoj111}) of Theorem \ref{karakterizacija-nenegativnih-kompleksnih-mer}:}
		Take $A,B\in (W_1)_+$.
		For $A,B$ with finite spectral decompositions the condition (\ref{pogoj111}) follows by (\ref{pogoj1111}). If not both $A, B$ have finite spectral 
		decompositions, then we have to prove that
			\bequ \label{enakost-101}
				\lim_{\ell\to\infty} \left(M_{S_\ell(A+B)}(\Delta)-M_{S_\ell(A)}(\Delta)-M_{S_\ell(B)}(\Delta)\right)=0,
			\eequ
		where the limit is taken in the norm topology and the sequence $S_{\ell}(A)$ (resp.~$S_{\ell}(B)$, $S_{\ell}(A+B)$) is a constant sequence if $A$ (resp.~$B$, $A+B$) has a finite spectral 
		decomposition, i.e., $S_{\ell}(A)=A$ for every $\ell\in\NN$.
		Define 
			$$T_{\ell}:=S_\ell(A+B)-S_\ell(A)-S_\ell(B)$$ 
		and notice $\left\|T_{\ell}\right\|\leq \frac{3}{\ell}.$
		Further on, for every $\epsilon>0$ there exists $N\in\NN$ such that for $j\geq N$, 
			$$\left\|T_{\ell}-S_j(T_{\ell})\right\|\leq \epsilon,\quad  
			\left\|M_{T_{\ell}}(\Delta)-M_{S_j(T_{\ell})}(\Delta)\right\|\leq \epsilon,$$
		where $S_j(T_{\ell})$ is defined as in (\ref{Riemannove-vsote-operatorja}).
		Therefore also 
			$$\left\|S_j(T_\ell)\right\|=\left\|T_\ell-T_\ell+S_j(T_\ell)\right\|\leq \left\|T_\ell\right\| + \left\|T_\ell-S_j(T_\ell)\right\|
			\leq \frac{3}{\ell}+\epsilon.$$
		As for the existence and uniqueness we estimate
			$\left\|M_{S_j(T_{\ell})}(\Delta)\right\|\leq \left(\frac{3}{\ell}+\epsilon\right)k_\Delta$ 
		and hence
			$\left\|M_{T_{\ell}}(\Delta)\right\| \leq \epsilon + \left(\frac{3}{\ell}+\epsilon\right)k_\Delta$.
		Since $\epsilon>0$ was arbitrary, we conclude that $$\left\|M_{T_{\ell}}(\Delta)\right\| \leq \frac{3\cdot k_{\Delta}}{\ell}.$$
		Hence, $\lim_{\ell} M_{T_{\ell}}(\Delta)=0$ which proves (\ref{enakost-101}).
		
	\textsl{Step 4 - $\{M_A\}_{A\in (W_{1})_+}$ satisfies the condition (\ref{pogoj222}) of Theorem \ref{karakterizacija-nenegativnih-kompleksnih-mer}:} Analoguous to the proof of the Step 3.	
			
	\textsl{Step 5 - $\{M_A\}_{A\in (W_{1})_+}$ satisfies the condition (\ref{pogoj333}) of Theorem \ref{karakterizacija-nenegativnih-kompleksnih-mer}:} Take $A\in (W_1)_+$.
				Each $S_{\ell}(A)=:\sum_i \lambda_{i\ell} P_{i\ell}$ has a finite spectral decomposition, where $\lambda_{i\ell}$ are real numbers and
				$P_{i\ell}$ mutually orthogonal hermitian projections.
				By $\left\|A-S_\ell(A)\right\|\leq \frac{1}{\ell}$ and $P_{i\ell} P_{j\ell}=0$ for $i\neq j$,
				we have 
					$\left\|A\right\|-\frac{1}{\ell}\leq \left\|\sum_{i=1}^{k_\ell} \lambda_{i\ell} P_{i\ell}\right\|\leq 
					\left\|A\right\|+\frac{1}{\ell}$
				and hence
					$\max_{i} \left|\lambda_{i\ell} \right| \in \left(\left\|A\right\|-\frac{1}{\ell},\left\|A\right\|+\frac{1}{\ell}\right).$
				It follows that
					\beqn
						&&\left\|M_A(\Delta)\right\| = \left\|\lim_{\ell\to\infty} M_{S_\ell(A)}(\Delta)\right\|
						=	\lim_{\ell\to\infty} \left\|\sum_{i=1}^{k_\ell}	\lambda_{i\ell} M_{P_{i\ell}}(\Delta)\right\|\\
						&\underbrace{\leq}_{(\ast)}&
																	\lim_{\ell\to\infty}\left(\max_{i=1} {\left|\lambda_{i\ell}\right|}
																	\left\|\sum_{i=1}^{k_\ell} M_{P_{i\ell}}(\Delta)\right\|\right)
																\leq \lim_{\ell\to\infty} \left(\left(\|A \|+\frac{1}{\ell} \right)
																\left\|\sum_{i=1}^{k_\ell} M_{P_{i\ell}}(\Delta)\right\|\right)	\\
																&\underbrace{\leq}_{ \text{by}\;(\ref{pogoj1111})}& 
																\lim_{\ell\to\infty}\left(\left(\left\|A\right\|+\frac{1}{\ell} \right)
																\left\|M_{\sum_{i=1}^k P_{i\ell}}(\Delta)\right\|\right)
																\underbrace{\leq}_{\substack{\text{by}\;(\ref{pogoj3333})}}
																\lim_{\ell\to\infty} \left(\left\|A\right\|+\frac{1}{\ell} \right)\cdot k_\Delta\\
																&=& \|A \|k_\Delta,
					\eeqn 
				where $(\ast)$ follows by $M_{P_{i\ell}}(\Delta)\succeq 0$ for every $i,\ell\in \NN$ and every $\Delta\in\cS$.
\eprf

As a corollary we obtain the following equivalent definiton of a non-negative spectral measure.

\bcor \label{prvotna-definicija}
	Let $(X,\cS,W_1,W_2)$ be a measure space. A set function $M \colon \cS \to B(W_1,W_2)$
	is a non-negative spectral measure if for every hermitian projection $P\in (W_1)_p$ the set functions $M_P$ 
	are spectral measures and the equality 
		$$M_P(\Delta_1)M_Q(\Delta_2)=M_{PQ}(\Delta_1\cap \Delta_2)$$
	holds for all hermitian projections $P,Q\in (W_1)_p$ and all sets $\Delta_1, \Delta_2\in \cS$.
\ecor

\brem
	Notice that by the definition, a non-negative spectral measure is a \emph{a non-negative measure} with the properties from Corollary \ref{prvotna-definicija}
	(see Section \ref{nenegativne-spektralne}).
	Corollary \ref{prvotna-definicija} shows that \emph{a set function} with these properties is automatically
	a non-negative measure. Hence, we obtained precisely the definition of a non-negative spectral measure from the Introduction.
\erem

\section{Integral representations of representations $\rho: C(X,W_1) \to W_2$}
\label{integralske-reprezentacije}

Let $(X,\B(X),W_1\subseteq B(\cH),W_2\subseteq B(\cK),M)$ be a space with a non-negative spectral measure $M$ (see Section \ref{nenegativne-spektralne}), where $X$ is a compact Hausdorff space and $\B(X)$ is a Borel $\sigma$-algebra on $X$. 
We call $M$ \textsl{regular} if the spectral measures $M_P$ are regular for every $P\in (W_1)_p$, i.e., complex measures
$$(M_P)_{k_1,k_2}:\B(X)\to \CC,\quad (M_P)_{k_1,k_2}(\Delta):=\left\langle M_P(\Delta)k_1,k_2\right\rangle$$ 
are regular for every $k_1, k_2\in \cK$ and every $P\in (W_1)_p$. 
$M$ is \textsl{normalized} if $M(\id_{\cH})=\id_{\cK}$, where $\id_\cH$, $\id_\cK$ denote the identity operators on $\cH$, $\cK$ respectively.

The main result of this article is the following.

\bthe \label{reprezentacijski-izrek-1}
	Let $X$, $W_1$, $W_2$ be as above and 
		$$\rho: C(X,\CC)\otimes W_1 \to W_2$$ 
	a bounded linear map.
	The following statements are equivalent.
	\benu
		\item $\rho: C(X,\CC)\otimes W_1 \to W_2$ is a unital algebra homomorphism such
			that $\rho(F^\ast)=\rho(F)^{\ast}$ for every $F\in C(X,\CC)\otimes W_1$.
		\item	There exists a unique regular normalized non-negative spectral measure $M:\B(X)\to B(W_1,W_2)$ such that
				\begin{equation} \label{def-rho} \rho(F)=\int_X F\; dM \end{equation}
			for every $F\in C(X,\CC)\otimes W_1$.
	\eenu
\ethe

\bprf
	Direction $(2)\Rightarrow (1).$ For $\rho$ satisfying (\ref{def-rho}) we have to prove the linearity, the multiplicativity of $\rho$ and the equality
	$\rho(F^\ast)=\rho(F)^{\ast}$ for every $F\in C(X,\CC)\otimes W_1$. The linearity follows by Proposition \ref{lastnosti-kompleksnega-integrala}, while
	the multiplicativity by Proposition \ref{integracija-po-multiplikativni-meri-je-mutliplikativna}. To show 
		$\int_X F^\ast\; dM=\left(\int_X 	F\;dM\right)^{\ast}$  
	it suffices, by the linearity, to consider elements of the form $F=f\otimes A$, $f\in \cI(M)_+$, $A\in (W_1)_+$.
	Since $M_A$ is a positive operator-valued measure, we have
	$\int_X (f\otimes A)^\ast\; dM=\int_X (f\otimes A)\;dM=\left(\int_X (f\otimes A)\;dM\right)^{\ast}$ and the result follows.
	
	Direction $(1)\Rightarrow (2).$ Since $\rho$ is an algebra homomorphism such that 
	$\rho(F^\ast)=\rho(F)^{\ast}$ for every $F\in C(X,\CC)\otimes W_1$, the maps
		$\rho_P: C(X,\CC) \to W_2,$ $\rho_P(f):=\rho(f\otimes P)$ are $\ast$-representations 
	for every $P\in (W_1)_p$.	
	By Theorem \ref{reprezentacija}, there exist unique spectral measures
	$F_P:\B(X)\to W_2$ such that $\rho_P(f)=\int_X f\; dF_P$ holds for every $f\in C(X,\CC)$ and every $P\in (W_1)_p$. The idea is to show that the family
	$\{F_P\}_{P\in (W_1)_p}$ satisfies the conditions of Theorem \ref{karakterizacija-nenegativnih-spektralnih-mer} 
	to obtain a non-negative spectral measure $M$ representing $\rho$.
	
	\textsl{The family $\{F_P\}_{P\in (W_1)_p}$ satisfies the condition (\ref{pogoj1111}) of Theorem \ref{karakterizacija-nenegativnih-spektralnih-mer}:}
	Let $P_i,Q_j\in (W_1)_p$ be hermitian projections and $\lambda_i,\mu_j \in\RR$ real numbers, such that $\sum_{i=1}^{n}\lambda_i P_i=\sum_{j=1}^m \mu_j Q_j$.
	We have to show that for every set $\Delta\in\B(X)$, the equality $$\sum_{i=1}^{n}\lambda_i F_{P_i}(\Delta)=\sum_{j=1}^m \mu_j F_{Q_j}(\Delta)$$ holds.
	Since the function $\chi_{\Delta}$ is a bounded Borel function,
	by Lemma \ref{density}, there is a net $\{f_k\}\subset C(X,\CC)$ such that $\int_X f_k\;d\mu \to \mu(\Delta)$ for every measure $\mu\in C(X,\CC)^{\ast\ast}$.
	Therefore for all $k_1, k_2\in\cK$
		\begin{eqnarray*}
			\left\langle \rho\left(f_k\otimes \sum_{i=1}^{n}\lambda_i P_i\right) k_1,k_2\right\rangle &=& 
			\left(\sum_{i=1}^{n}\lambda_i \int_X f_k \; d(F_{P_i})_{k_1,k_2}\right) \to\\ 
				&\to&	\sum_{i=1}^{n}\lambda_i \int_X \chi_\Delta \; d(F_{P_i})_{k_1,k_2} = \sum_{i=1}^{n}\lambda_i (F_{P_i})_{k_1,k_2}(\Delta).
		\end{eqnarray*}
	and analoguously 
	$\left\langle \rho\left(f_k\otimes \sum_{j=1}^{m}\mu_j Q_j\right) k_1,k_2\right\rangle \to \sum_{j=1}^{m}\mu_j(F_{Q_j})_{k_1,k_2}(\Delta)$.
	Since $\rho\left(f_k\otimes \sum_{i=1}^{n}\lambda_i P_i\right)=\rho\left(f_k\otimes \sum_{j=1}^{m}\mu_j Q_j\right)$
	holds for every $k\in \NN$, 
	it follows that $\sum_{i=1}^{n}\lambda_i F_{P_i}(\Delta) = \sum_{j=1}^{m}\mu_j F_{Q_j}(\Delta).$
	
	\textsl{The family $\{F_P\}_{P\in (W_1)_p}$ satisfies the condition (\ref{pogoj3333}) of Theorem \ref{karakterizacija-nenegativnih-spektralnih-mer}:}
	Let $P\in (W_1)_p$ be a hermitian projection and $\Delta\in\Bor(X)$ a Borel set. We have to find a constant $k_\Delta\in\RR^{>0}$ such that  
	$\left\|F_P(\Delta)\right\|\leq k_\Delta$. We know that 
		\beqn
			\|F_P(X)\| &=& \left\|\int_X 1\; dF_P\right\|  = \left\|\rho_P(1)\right\|
				= \left\|\rho(1\otimes P)\right\| \leq \left\|\rho\right\|\left\|1\otimes P\right\|_{\infty}
												          = \left\|\rho\right\|,
		\eeqn
	where we used the continuity of $\rho$ for the inequality. 
	By the finite additivity of $F_P$, it follows that $\left\|F_P(\Delta)\right\| \leq \left\|\rho\right\|$ for every $\Delta\in\Bor(X)$.
	
	\textsl{The family $\{F_P\}_{P\in (W_1)_p}$ satisfies the condition (\ref{pogoj4444}) of Theorem \ref{karakterizacija-nenegativnih-spektralnih-mer}:}
	Let $P,Q\in (W_1)_p$ be hermitian projections and $\Delta_1,\Delta_2\in \Bor(X)$ Borel sets. We have to show that 
	$$M_P(\Delta_1)M_Q(\Delta_2)=M_{PQ}(\Delta_1\cap \Delta_2).$$
	By Lemma \ref{density}, there exists a net $\{f_k\}\subset C(X,\CC)$, such that $\int_X f_k \cdot g\; d\mu \to \int_X \chi_{\Delta_1} \cdot g\; d\mu$ 
	for every $\mu\in C(X,\CC)^{\ast\ast}$ and every bounded Borel function $g$. 
	Therefore $\int_X f_k \cdot g\; dM_{PQ} \to  \int_X \chi_{\Delta_1} \cdot g\; dM_{PQ}$ 
	and $\int_X f_k \; dM_{P} \to  \int_X \chi_{\Delta_1} \; dM_{P}$ in the weak operator topology.
	Hence for $g\in C(X,\CC)$ 
		\beqn 
		\int_X \chi_{\Delta_1}\cdot g\; dM_{PQ} &=& \lim \int_X f_k \cdot g\; dM_{PQ} = \lim \rho(f_k\cdot g\otimes PQ)\\
			&=& \lim \rho(f_k\otimes P) \rho(g\otimes Q)
			=   \lim \left(\int_X f_k \; dM_{P}\right) \left(\int_X g \; dM_{Q}\right)\\
			&=& \left(\int_X \chi_{\Delta_1}\;dM_{P}\right) \left(\int_X  g\; dM_{Q}\right),
		\eeqn
	where all the limits are in the weak operator topology.
	By Lemma \ref{density}, there exists a net $\{g_k\}\subset C(X,\CC)$, such that 
	$\int_X \chi_{\Delta_1} \cdot g_k\; d\mu \to \int_X \chi_{\Delta_1} \cdot \chi_{\Delta_2}\; d\mu$
	and $\int_X g_k\; d\mu \to \int_X \chi_{\Delta_2}\;d\mu$ for every $\mu\in C(X,\CC)^{\ast\ast}$.
	Therefore
		\beqn 
			\int_X \chi_{\Delta_1} \cdot \chi_{\Delta_2}\; dM_{PQ} &=& \lim \int_X \chi_{\Delta_1}\cdot g_k\; dM_{PQ} \\
			&=& \lim \left(\int_X \chi_{\Delta_1} \; dM_{P}\right) \left(\int_X g_k \; dM_{Q}\right)\\
			&=& \left(\int_X \chi_{\Delta_1} \; dM_{P}\right) \left(\int_X \chi_{\Delta_2} \; dM_{Q}\right),
		\eeqn
	where all the limits are in the weak operator topology.
	It follows that
		$M_{PQ}(\Delta_1\cap \Delta_2)=M_P(\Delta_1)M_Q(\Delta_2).$
	
	\textsl{$M$ is the representing measure of $\rho$:}
	By the linearity and the continuity of $\rho$ and $\int$, it suffices to consider the elements $F\in C(X,\CC)\otimes (W_1)_p$ of the form $f\otimes P$.	By the construction of the measures $F_P$, we have
		$\rho(f\otimes P)=\int_Xf\;dF_P=\int_X(f\otimes P)\;dM$. Hence $M$ represents $\rho$.
	
	\textsl{$M$ is unique, regular and normalized:} This follows from the uniqueness and the regularity of each $F_P$ and the unitality of
	$\rho$.
\eprf

Corollary \ref{reprezentacijski-izrek-2} is a slight generalization of Theorem \ref{reprezentacija}, i.e., the map $\rho: C(X,\CC)\otimes W_1 \to W_2$
is replaced by the map $\rho: C(X,W_1)\to W_2$. Note that we also do not need the boundedness of $\rho$ in the statement of Corollary \ref{reprezentacijski-izrek-2}, since it automatically follows by (\ref{prva-tocka}) or (\ref{druga-tocka}).

\bcor \label{reprezentacijski-izrek-2}
	Let $X$ be a compact Hausdorff space, $\B(X)$ a Borel $\sigma$-algebra on $X$, $W_1$, $W_2$ von Neumann algebras and 
		$$\rho: C(X,W_1) \to W_2$$ a linear map.
	The following statements are equivalent.
	\benu
		\item \label{prva-tocka} $\rho: C(X,W_1) \to W_2$ is a unital $\ast$-representation. 
		\item	\label{druga-tocka} There exists a unique regular normalized non-negative spectral measure $M:\B(X)\to B(W_1,W_2)$ such that
				$$\rho(F)=\int_X F\; dM$$
			for every $F\in C(X,W_1)$.
	\eenu
\ecor

\bprf
	Direction $(2)\Rightarrow (1)$: We have to proof that $\rho$ is multiplicative and satisfies the $\ast$-condiditon. 
	For the elements from the set $C(X,\CC)\otimes W_1$ the proof is the same as the proof of $(2)\Rightarrow (1)$
	in Theorem \ref{reprezentacijski-izrek-1}. By the definition of the integration with respect to $M$ on a general element from the set $C(X,W_1)$ (see 
	(\ref{razsirimo-integracijo})), it also follows for all elements from $C(X,W_1)$.
	
	For the direction $(1)\Rightarrow (2)$ we first notice that $\rho$ is bounded by \cite[4.8.~Theorem., p.~247]{Con}. 
	Then we apply Theorem 
	\ref{reprezentacijski-izrek-1} to $\rho|_{C(X,\CC)\otimes W_1}$ to obtain the unique representing measure $M$ for 
	$\rho|_{C(X,\CC)\otimes W_1}$. 
	For a general element $F$ from $C(X,W_1)$ there exists a sequence 
	$F_i\subseteq C(X,\CC)\otimes W_1$ such that $\lim_i F_i=F$, where the limit is taken in the supremum norm.
	By the continuity of $L$, we have $L(F)=\lim_i L(F_i)$, where the limit is taken in the usual operator norm.
	Again by the definition of the integration with respect to $M$ on a general element from the set $C(X,W_1)$ (see (\ref{razsirimo-integracijo})), $\int_X F\;dM=\lim_{i} \int_X F_i\;dM$.
	Hence $L(F)=\int_X F\;dM$.
\eprf

\brem
	As in Remark \ref{regularnost}.(1), there is a corresponding version of Corollary \ref{reprezentacijski-izrek-2}
	for a locally compact Hausdorff space $X$ and linear maps of the form $\rho: C_0(X,W_1) \to W_2$, 
	where the set $C_0(X,W_1)$ denotes elements of $C(X,W_1)$ vanishing at infinity, i.e., $F\in C_0(X,W_1)$ iff for every $\epsilon>0$ there exists a compact set $K_\epsilon$, such
	that $\left\|F(x)\right\|<\epsilon$ for every $x\in K_\epsilon^c$.
	Namely, for a locally compact Hausdorff space $X$, unital $\ast$-representations
	$\rho: C_0(X,W_1) \to W_2$ are in one-to-one correspondence with the regular normalized non-negative
	spectral measures $M:\Bor(X_{\infty}) \to B(W_1,W_2)$, where $X_{\infty}$ stands for the one point compactification of $X$.
\erem

\noindent \textbf{Acknowledgement.} I would like to thank to my advisor Jaka Cimprič for many useful discussions
while studying the material above and suggesting improvements of the manuscript. I am also grateful to Bojan
Magajna for the discussions concerning von Neumann algebras, vector measures and Remark \ref{regularnost}.(1).

\end{document}